\input amstex
\documentstyle{amsppt}
\magnification=\magstep1
 \hsize 13cm \vsize 18.35cm \pageno=1
\loadbold \loadmsam
    \loadmsbm
    \UseAMSsymbols
\topmatter
\NoRunningHeads
\title The fermionic $p$-adic integrals on $\Bbb Z_p$ associated
with extended $q$-Euler numbers and polynomials
\endtitle
\author
  Taekyun Kim
\endauthor

\abstract The purpose of this paper is to present a systemic study
of some families of $q$-Euler numbers and polynomials of
N$\ddot{o}$rlund's type by using multivariate fermionic $p$-adic
integral on $\Bbb Z_p$. Moreover, the study of these higher-order
$q$-Euler numbers and polynomials implies some interesting
$q$-analogue of stirling numbers identities.

\endabstract
\thanks  2000 AMS Subject Classification: 11B68, 11S80
%\newline keywords and phrases
\newline  The present Research has been conducted by the research
Grant of Kwangwoon University in 2008
\endthanks
\endtopmatter

\document

{\bf\centerline {\S 1. Introduction/ Preliminaries}}

 \vskip 15pt
Let $p$ be a fixed odd prime number. Throughout this paper $\Bbb
Z_p$, $\Bbb Q_p$, $\Bbb C$ and $\Bbb C_p$ will, respectively, denote
the ring of $p$-adic rational integers, the field of $p$-adic
rational numbers, the complex number field and the completion of
algebraic closure of $\Bbb Q_p$. Let $v_p$ be the normalized
exponential valuation of $\Bbb C_p$ with
$|p|_p=p^{-v_p(p)}=\frac{1}{p}.$ When one talks of $q$-extension,
$q$ is variously considered  as and indeterminate, a complex number
$q\in \Bbb C$ or $p$-adic number $q\in \Bbb C_p$. If $q\in \Bbb C$,
one normally assumes $|q|<1$. If $q\in\Bbb C_p$, one normally
assumes $|1-q|_p<1.$ We use the notation
$$[x]_q=\frac{1-q^x}{1-q}, \text{ and }
[x]_{-q}=\frac{1-(-q)^x}{1+q}, \text{ (see [1-10])}.$$ The
$q$-factorial is defined as $[n]_q!=[n]_q[n-1]_q\cdots[2]_q[1]_q$
and the Gaussian binomial coefficient is also defined by
$${\binom{n}{k}}_q=\frac{[n]_q!}{[n-k]_q![k]_q!}=\frac{[n]_q [n-1]_q\cdots [n-k+1]_q}{[k]_q!}, \text{ (see [9])}.\tag1$$
Note that
$$\lim_{q\rightarrow
1}{\binom{n}{k}}_q=\binom{n}{k}=\frac{n!}{(n-k)!k!}=\frac{n(n-1)\cdots(n-k+1)}{k!}.$$
From (1), we easily note that
$${\binom{n+1}{k}}_q={\binom{n}{k-1}}_q+q^k{\binom{n}{k}}_q=q^{n-k}{\binom{n}{k-1}}_q+{\binom{n}{k}}_q, \text{ (see [5])}.$$
The $q$-binomial formulae are known that

$$(b;q)_n=(1-b)(1-bq)\cdots(1-bq^{n-1})=\sum_{i=0}^n
{\binom{n}{i}}_q q^{\binom{i}{2}}(-1)^ib^i, \tag2$$ and

$$\frac{1}{(b;q)_n}=\frac{1}{(1-b)(1-bq)\cdots(1-bq^{n-1})}=\sum_{i=0}^{\infty}{\binom{n+i-1}{i}}_q
b^i.$$ We say that $f$ is uniformly differentiable function at a
point $a\in\Bbb Z_p$, and write $f\in UD(\Bbb Z_p),$ if the
difference quotient $F_{f}(x, y)=\frac{f(x)-f(y)}{x-y}$ have a limit
$f^{\prime}(a)$ as $(x,y)\rightarrow (a,a).$ For $f \in UD(\Bbb
Z_p),$ the fermionic $p$-adic $q$-integral on $\Bbb Z_p$ is defined
as
$$I_{-q}(f)=\int_{\Bbb Z_p}f(x) d\mu_{-q}(x)=\lim_{N\rightarrow
\infty}\frac{1+q}{1+q^{p^N}}\sum_{x=0}^{p^N-1}f(x)(-q)^x, \text{
(see [9])}. \tag3$$ From (3),  we derive the fermionic  $p$-adic
integral on $\Bbb Z_p$  as follows.
$$\lim_{q\rightarrow 1}I_{-q}(f)=I_{-1}(f)=\int_{\Bbb Z_p} f(x)
d\mu_{-1}(x). \tag4$$ For $n\in \Bbb N$, let $f_n(x)= f(x+n)$, we
have
$$I_{-1}(f_n)=(-1)^nI_{-1}(f)+2\sum_{l=0}^{n-1}(-1)^{n-1-l}f(l). \tag5$$
From (5), we can easily derive the Euler polynomials, $E_n(x),$ as
follows.
$$\int_{\Bbb Z_p} e^{(x+y)t}d\mu_{-1}d\mu_{-1}(y)=\frac{2}{e^t
+1}e^{xt}=\sum_{n=0}^{\infty}E_n(x)\frac{t^n}{n!}, \text{ (see
[1-23])}.\tag6$$ Note that $E_n(0)=E_n$ are called the $n$-th Euler
numbers. Now, we consider the Euler polynomials of
N$\ddot{o}$rlund's type as follows.
$$\int_{\Bbb Z_p}\cdots\int_{\Bbb Z_p}e^{(x+x_1+\cdots+x_r)t}d\mu_{-1}(x_1)\cdots d\mu_{-1}(x_r)=\left(\frac{2}{e^t+1} \right)^r e^{xt}=\sum_{n=0}^{\infty}
E_n^{(r)}(x) \frac{t^n}{n!},\tag7$$ and
$$\left( \frac{e^t+1}{2} \right)^r
e^{xt}=\sum_{n=0}^{\infty}E_n^{(-r)}(x)\frac{t^n}{n!}, \text{ (see
[9])}. \tag 7-1$$ In the special case $x=0$,
$E_n^{(-r)}(0)=E_n^{(r)}$ and $E_n^{(r)}(0)=E^{(-r)}$ are called the
Euler numbers of N$\ddot{o}$rlund's type. Let $(Eh)(x)=h(x+1)$ be
the shift operator. Then the $q$-difference operator $\Delta_q$ is
defined as
$$\Delta_q^n=\prod_{i=1}^n(E-q^{i-1}I), \text{ where
$(Ih)(x)=h(x)$, (see [5, 10]).} \tag8$$ From (8), we note that
$$f(x)=\sum_{n\geq 0}{\binom{x}{n}}_q \Delta_{q}^nf(0), \tag9$$
where
$$\Delta_q^n
f(0)=\sum_{k=0}^n{\binom{n}{k}}_q(-1)^kq^{\binom{k}{2}}f(n-k),
\text{( see [5, 10])}.$$ The $q$-stirling number of the second kind
(as defined by Carlitz) is given by
$$S_2(n,k;q)=\frac{q^{-\binom{k}{2}}}{[k]_q!}\sum_{j=0}^k(-1)^jq^{\binom{j}{2}}{\binom{k}{j}}_q[k-j]_q^n,
\text{ (see [10])}.\tag10$$ By (9) and (10), we see that
$$S_2(n,k;q)=\frac{q^{-\binom{k}{2}}}{[k]_q!}\Delta_q^k0^n, \text{
(see [10])}.\tag11$$ In this paper, the $q$-extension of (7) are
variously considered. From these $q$-extensions, we derive some
interesting identities and relations of Euler polynomials and
numbers of N$\ddot{o}$rlund's type. The purpose of this paper is to
present  a systemic study of some families of $q$-Euler numbers and
polynomials of N$\ddot{o}$rlund's type by using  multivariate
fermionic $p$-adic integral on $\Bbb Z_p$

\vskip 10pt

{\bf\centerline {\S 2. The $q$-extension of Euler numbers and
polynomials  of  $N\ddot{o}rlund$ type}} \vskip 10pt

In this section, we assume that $q\in \Bbb C_p$ with $|1-q|_p<1$.
First, we consider the $q$-extension of (6) as follows.
$$\aligned
\sum_{n=0}^{\infty}E_{n,q}(x)\frac{t^n}{n!}&=\int_{\Bbb
Z_p}e^{[x+y]_qt}d\mu_{-q}(x)=\sum_{n=0}^{\infty}\frac{2}{(1-q)^n}\sum_{l=0}^{n}\left(\frac{\binom{n}{l}(-1)^lq^{lx}}{1+q^l}\right)
\frac{t^n}{n!}\\
&=2\sum_{m=0}^{\infty}(-1)^me^{[m+x]_qt}.\endaligned$$ Therefore, we
obtain the following lemma.

 \proclaim{ Lemma 1}
For $ n\geq 0$,  we have
$$E_{n, q}(x)=2\sum_{m=0}^{\infty}(-1)^m
[m+x]_q^n=\frac{2}{(1-q)^n}\sum_{l=0}^{n}\frac{\binom{n}{l}(-1)^lq^{lx}}{1+q^l}.\tag12$$
\endproclaim
From (11), we note that

$$\aligned
[x]_q^n&=\sum_{k=0}^n{\binom{n}{k}}_q[k]_q!S_2(k,
n-k;q)q^{\binom{k}{2}}\\
& =\sum_{k=0}^n[x]_q[x-1]_q\cdots[x-k+1]_q
\frac{q^{\binom{k}{2}-\binom{n-k}{2}}}{[n-k]_q!}\Delta_q^{n-k}0^k\\
&=\sum_{k=0}^n\frac{q^{\binom{k}{2}-\binom{n-k}{2}}}{[n-k]_q!}\Delta_q^{n-k}0^k\frac{1}{(1-q)^k}
\sum_{l=0}^k {\binom{k}{l}}_qq^{\binom{l}{2}}(-1)^lq^{l(x-k+1)}.\\
\endaligned$$
Thus, we have
$$E_{n,q}=\sum_{k=0}^n\frac{q^{\binom{k}{2}}S_2(k,n-k;q)}{(1-q)^k}\sum_{l=0}^k{\binom{k}{l}}_qq^{\binom{l}{2}}(-1)^l\sum_{m=0}^l
\binom{l}{m}(q-1)^mE_{m,q}(1-k).$$

Therefore, we obtain the following theorem.

\proclaim{ Theorem 2} For $n\geq 0$, we have
$$E_{n,q}=\sum_{k=0}^n\frac{q^{\binom{k}{2}}S_2(k,n-k;q)}{(1-q)^k}\sum_{l=0}^k{\binom{k}{l}}_qq^{\binom{l}{2}}(-1)^l\sum_{m=0}^l
\binom{l}{m}(q-1)^mE_{m,q}(1-k).$$
\endproclaim
Let us consider the $q$-extension of  Eq.(7) as follows.
$$\aligned
E_{n,q}^{(r)}(x)&=\int_{\Bbb Z_p}\cdots\int_{\Bbb
Z_p}[x+x_1+\cdots+x_r]_q^nd\mu_{-1}(x_1)\cdots d\mu_{-1}(x_r)\\
&=\frac{2^r}{(1-q)^n}\sum_{l=0}^n\binom{n}{l}(-1)^lq^{lx}\left(\frac{1}{1+q^l}\right)^r\\
&=2^r\sum_{m=0}^{\infty}\binom{m+r-1}{m}(-1)^m[m+x]_q^n.
\endaligned\tag13$$
Let
$F_q^{(r)}(t,x)=\sum_{n=0}^{\infty}E_{n,q}^{(r)}(x)\frac{t^n}{n!}.$
Then we have
$$F_q^{(r)}(t,x)=2^r\sum_{m=0}^{\infty}\binom{m+r-1}{m}(-1)^me^{[m+x]_qt}.\tag14$$
In the special case $x=0$, $E_{n,q}^{(r)}(0)=E_{n,q}^{(r)}$ are
called the $q$-extension of Euler numbers of order $r$. In the sense
of the $q$-extension of Eq.(7-1), we consider the $q$-extension of
Euler polynomials of N$\ddot{o}$rlund's type as follows.
$$G_q^{(r)}(t,x)=F_q^{(-r)}(t,x)=\frac{1}{2^r}\sum_{m=0}^r\binom{r}{m}e^{[m+x]_qt}=\sum_{n=0}^{\infty}E_{n,q}^{(-r)}(x)\frac{t^n}{n!}.\tag15$$
By (15), we see that
$$E_{n,q}^{(-r)}(x)=\frac{1}{2^r}\sum_{m=0}^{r}\binom{r}{m}[m+x]_q^n.$$
Therefore, we obtain the following theorem.

\proclaim{ Theorem 3} For $r\in \Bbb N$, $n\geq 0$,  let
$$2^r\sum_{m=0}^{\infty}\binom{m+r-1}{m}(-1)^me^{[m+x]_qt}=\sum_{n=0}^{\infty}E_{n,q}^{(r)}(x)\frac{t^n}{n!}.$$
Then we have
$$\aligned
E_{n,q}^{(r)}(x)&=\frac{2^r}{(1-q)^n}\sum_{l=0}^n\binom{n}{l}(-1)^lq^{lx}\left(\frac{1}{1+q^l}\right)^r\\
&=2^r\sum_{m=0}^{\infty}\binom{m+r-1}{m}(-1)^m[m+x]_q^n,
\endaligned$$
and
$$\aligned
E_{n,q}^{(-r)}(x)&=\frac{1}{2^r(1-q)^n}\sum_{l=0}^n\binom{n}{l}(-1)^lq^{lx}\left(1+q^l\right)^r\\
&=\frac{1}{2^r}\sum_{m=0}^{r}\binom{r}{m}[m+x]_q^n.
\endaligned$$
\endproclaim
$E_{n,q}^{(-r)}(0)=E_{n,q}^{(-r)}$ are called the $q$-extension of
Euler numbers of N$\ddot{o}$rlund's type. For $h\in\Bbb Z$,
$r\in\Bbb N$, let us define the extended higher-order $q$-Euler
polynomials as follows.
$$E_{n,q}^{(h,r)}(x)=\int_{\Bbb Z_p}\cdots \int_{\Bbb
Z_p}q^{\sum_{j=1}^r(h-j)x_j}[x+x_1+\cdots+x_r]_q^nd\mu_{-1}(x_1)\cdots
d\mu_{-1}(x_r). \tag16$$ Then, we have
$$\aligned
E_{n,q}^{(h,r)}(x)&=\frac{2^r}{(1-q)^n}\sum_{l=0}^{n}\frac{\binom{n}{l}(-1)^lq^{lx}}{(-q^{h-1+l};q^{-1})_r}
=\frac{2^r}{(1-q)^n}\sum_{l=0}^n\frac{\binom{n}{l}(-1)^lq^{lx}}{(-q^{h-r+l};q)_r}\\
&=2^r\sum_{m=0}^{\infty}{\binom{m+r-1}{m}}_qq^{(h-r)m}(-1)^m[x+m]_q^n.
\endaligned\tag17$$
Let
$F_q^{(h,r)}(t,x)=\sum_{m=0}^{\infty}E_{n,q}^{(h,r)}(x)\frac{t^n}{n!}.$
Then we easily see that
$$F_q^{(h,r)}(t,x)=2^r\sum_{m=0}^{\infty}q^{(h-r)m}(-1)^m
{\binom{m+r-1}{m}}_q e^{[m+x]_qt}. \tag18$$

Therefore we obtain the following theorem.

\proclaim{Theorem 4} For $h\in\Bbb Z, n\geq 0$, let
$$2^r\sum_{m=0}^{\infty}
q^{(h-r)m}(-1)^m{\binom{m+r-1}{m}}_qe^{[m+x]_qt}=\sum_{n=0}^{\infty}E_{n,q}^{(h,r)}(x)\frac{t^n}{n!}.$$
Then we have
$$\aligned
E_{n,q}^{(h,r)}(x)&=\frac{2^r}{(1-q)^n}\sum_{l=0}^n\frac{\binom{n}{l}(-1)^lq^{lx}}{(-q^{h-r+l};q)_r}\\
&=2^r\sum_{m=0}^{\infty}q^{(h-r)m}(-1)^m{\binom{m+r-1}{m}}_q
[x+m]_q^n.
\endaligned$$
\endproclaim
Now, we define the extended higher-order N$\ddot{o}$rlund's type
$q$-Euler polynomials as follows.

$$E_{n,q}^{(h,-r)}(x)=\frac{1}{(1-q)^n}\sum_{l=0}^n\frac{\binom{n}{l}(-1)^lq^{lx}}{\int_{\Bbb
Z_p}\cdots \int_{\Bbb
Z_p}q^{l(x_1+\cdots+x_r)}q^{\sum_{j=1}^r(h-j)x_j}d\mu_{-1}(x_1)\cdots
d\mu_{-1}(x_r)}. \tag19$$

 From (19), we note that
 $$\aligned
E_{n,q}^{(h,-r)}(x)&=\frac{1}{2^r(1-q)^n}\sum_{l=0}^n\binom{n}{l}(-1)^lq^{lx}(-q^{h-r+l};q)_r\\
&=\frac{1}{2^r}\sum_{m=0}^r {\binom{r}{m}}_q
q^{\binom{m}{2}}q^{(h-r)m}[m+x]_q^n.
 \endaligned\tag20$$

Let $F_q^{(h,-r)}(t,x)=\sum_{n=0}^{\infty}E_{n,q}^{(h,-r)}(x)
\frac{t^n}{n!}.$ Then we have
$$F_q^{(h,-r)}(t,x)=\frac{1}{2^r}\sum_{m=0}^rq^{\binom{m}{2}}q^{(h-r)m}
{\binom{r}{m}}_q e^{[m+x]_qt}. \tag21$$

Therefore, we obtain the following theorem.

\proclaim{ Theorem 5} For $h\in \Bbb Z$, $n\geq 0$, $r\in\Bbb N$,
let
$$\frac{1}{2^r}\sum_{m=0}^{\infty}q^{\binom{m}{2}}q^{(h-r)m}{\binom{r}{m}}_q
e^{[m+x]_qt}=\sum_{n=0}^{\infty}E_{n,q}^{(h,-r)}(x)\frac{t^n}{n!}.$$
Then we have
$$\aligned
E_{n,q}^{(h,-r)}(x)&=\frac{1}{2^r(1-q)^n}\sum_{l=0}^n\binom{n}{l}(-1)^lq^{lx}(-q^{h-r+l};q)_r\\
&=\frac{1}{2^r}\sum_{m=0}^rq^{\binom{m}{2}}q^{(h-r)m}{\binom{r}{m}}_q[m+x]_q^n.
\endaligned$$
\endproclaim

For $h=r$, we have
$$E_{n,q}^{(r,r)}(x)=\frac{2^r}{(1-q)^n}\sum_{l=0}^{n}\frac{\binom{n}{l}(-1)^lq^{lx}}{(-q^l;q)_r}
=2^r\sum_{m=0}^{\infty}{\binom{m+r-1}{m}}_q  (-1)^m [x+m]_q^n
,\tag22$$ and
$$\aligned
E_{n,q}^{(r,-r)}(x)&=\frac{1}{2^r (1-q)^n}\sum_{l=0}^n
\binom{n}{l}(-1)^lq^{lx}(-q^l;q)_r\\
&=\frac{1}{2^r}\sum_{m=0}^r q^{\binom{m}{2}}{\binom{r}{m}}_q
[m+x]_q^n.
\endaligned\tag23$$
It is easy to se that
$$\aligned
&\frac{q^{mx}2^r}{(-q^{m-r};q)_r}=\int_{\Bbb Z_p}\cdots \int_{\Bbb
Z_p}q^{\sum_{j=1}^r(m-j)x_j+mx} d\mu_{-q}(x_1)\cdots
d\mu_{-q}(x_r)\\
&=\int_{\Bbb Z_p}\cdots \int_{\Bbb
Z_p}\left([x+x_1+\cdots+x_r]_q(q-1)+1\right)^mq^{-\sum_{j=1}^r jx_j}
d\mu_{-1}(x_1)\cdots d\mu_{-1}(x_r)\\
&=\sum_{l=0}^m \binom {m}{l}(q-1)^l\int_{\Bbb Z_p}\cdots \int_{\Bbb
Z_p}[x+x_1+\cdots+x_r]_q^l q^{-\sum_{j=1}^r j
x_j}d\mu_{-1}(x_1)\cdots d\mu_{-1}(x_r)\\
&=\sum_{l=0}^m\binom{m}{l}(q-1)^lE_{l,q}^{(0,r)}(x).
\endaligned\tag24$$
By (24), we see that
$$\frac{q^{mx}2^r}{(-q^{m-r};q)_r}=\sum_{l=0}^m\binom{m}{l}(q-1)^lE_{l,q}^{(0,r)}(x).$$
It is known that
$$I_{-1}(f_1)+I_{-1}(f)=2f(0), \text{ where $f_1(x)=f(x+1)$}.\tag25$$
From (25), we derive
$$\aligned
&q^{h-1}\int_{\Bbb Z_p}\cdots\int_{\Bbb
Z_p}[x+1+x_1+\cdots+x_r]_q^mq^{\sum_{j=1}^r(h-j)x_j}d\mu_{-1}(x_1)\cdots
d\mu_{x_r}\\
&=-\int_{\Bbb Z_p}\cdots\int_{\Bbb
Z_p}[x+x_1+\cdots+x_r]^nq^{\sum_{j=1}^r(h-j)x_j}d\mu_{-1}(x_1)\cdots
d\mu_{-1}(x_r)\\
&+2\int_{\Bbb Z_p}\cdots\int_{\Bbb
Z_p}[x+x_2+\cdots+x_r]_q^nq^{\sum_{j=1}^{r-1}(h-1-j)x_{j+1}}d\mu_{-1}(x_2)\cdots
d\mu_{-1}(x_r).
\endaligned\tag26$$
By (26), we see that
$$q^{(h-1)}E_{n,q}^{(h,r)}(x+1)+E_{n,q}^{(h,r)}(x)=2E_{n,q}^{(h-1,r-1)}(x).\tag27$$
By simple calculation, we see that

$$\aligned
&q^x\int_{\Bbb Z_p}\cdots\int_{\Bbb
Z_p}[x_1+\cdots+x_r+x]_q^nq^{\sum_{j=1}^r(h-j+1)x_j}d\mu_{-1}(x_1)\cdots
d\mu_{-1}(x_r)\\
&=(q-1)\int_{\Bbb Z_p}\cdots\int_{\Bbb
Z_p}[x_1+\cdots+x_r+x]_q^{n+1}q^{\sum_{j=1}^r(h-j)x_j}d\mu_{-1}(x_1)\cdots d\mu_{-1}(x_r)\\
&+\int_{\Bbb Z_p}\cdots \int_{\Bbb
Z_p}q^{\sum_{j=1}^r(h-j)x_j}[x_1+\cdots+x_r+x]_q^nd\mu_{-1}(x_1)\cdots
d\mu_{-1}(x_r).
\endaligned\tag28$$
By (28), we see that
$$q^xE_{n,q}^{(h+1,r)}(x)=(q-1)E_{n+1,
q}^{(h,r)}(x)+E_{n,q}^{(h,r)}(x).$$

Therefore, we obtain the following
 proposition.

\proclaim{Proposition 6} For $h\in\Bbb Z$, $r\in\Bbb N$, $n\geq 0$,
we have
$$q^{(h-1)}E_{n,q}^{(h,r)}(x+1)+E_{n,q}^{(h,r)}(x)=2E_{n,q}^{(h-1,r-1)}(x),
$$
and
$$q^xE_{n,q}^{(h+1,r)}(x)=(q-1)E_{n+1,q}^{(h,r)}(x)+E_{n,q}^{(h,
r)}(x).$$ Moreover,
$$\frac{q^{mx}
2^r}{(-q^{m-r};q)_r}=\sum_{l=0}^m\binom{m}{l}(q-1)^lE_{l,q}^{(0,r)}(x).$$
\endproclaim

From (22), we note that
$$\aligned
E_{n,
q^{-1}}^{(r,r)}(r-x)&=\frac{2^r}{(1-q^{-1})^n}\sum_{l=0}^n\frac{\binom{n}{l}(-1)^lq^{-l(r-x)}}{(-q^{-l};q^{-1})_r}\\
&=(-1)^nq^{n+\binom{r}{2}}\frac{2^r}{(1-q)^n}\sum_{l=0}^n
\frac{\binom{n}{l}(-1)^lq^{lx}}{(-q^l;q)_r}\\
&=(-1)^nq^{n+\binom{r}{2}}E_{n,q}^{(r,r)}(x).
\endaligned$$
Hence, we have
$$\aligned
&\int_{\Bbb Z_p}\cdots \int_{\Bbb
Z_p}[r-x+x_1+\cdots+x_r]_{q^{-1}}^nq^{-\sum_{j=1}^r(r-j)x_j}d\mu_{-1}(x_1)\cdots
d\mu_{-1}(x_r)\\
&=(-1)^nq^{n+\binom{r}{2}}\int_{\Bbb Z_p}\cdots\int_{\Bbb
Z_p}[x+x_1+\cdots+x_r]_q^nq^{\sum_{j=1}^r(r-j)x_j}d\mu_{-1}(x_1)\cdots
d\mu_{-1}(x_r).
\endaligned$$
For $h=r$, we see that
$$E_{n,q^{-1}}^{(r,r)}(0)=(-1)^nq^{n+\binom{r}{2}}E_{n,q}^{(r,r)}(k).$$

From (27), we can also derive
$$q^{r-1}E_{n,q}^{(r,r)}(x+1)+E_{n,q}^{(r,r)}(x)=[2]_qE_{n,q}^{(r-1,r-1)}(x).$$

The stirling numbers of the first kind are defined as
$$\prod_{k=1}^n(1+[k]_qz)=\sum_{k=0}^nS_1(n,k;q)z^k, \text{ (see
[10])},\tag29$$ and
$$q^{\binom{m}{2}}{\binom{r}{m}}_q=\frac{q^{\binom{m}{2}}[r]_q\cdots[r-m+1]_q}{[m]_q!}=\frac{1}{[m]_q!}
\prod_{k=0}^{m-1}([r]_q-[k]_q).\tag30$$

It is easy to check that
$$\aligned
\prod_{k=0}^{n-1}\left(z-[k]_q \right)&=z^n\prod_{k=0}^{n-1}\left(1-\frac{[k]_q}{z}\right)\\
&=\sum_{k=0}^nS_1(n-1, k;q)(-1)^kz^{n-k}.
\endaligned\tag31$$
By (30) and (31), we see that
$$\prod_{k=0}^{m-1}\left([r]_q-[k]_q\right) =\sum_{k=0}^m S_1(m-1, k;q)(-1)^k[r]_q^{m-k}.\tag32$$
By(23) and (32), we obtain the following theorem.

\proclaim{ Proposition 7} For $r\in\Bbb N$, $n\in \Bbb Z_{+}$, we
have
$$E_{n,q}^{(r, -r)}(x)=\frac{1}{2^r [m]_q!}\sum_{m=0}^r\sum_{k=0}^mS_{1}(m-1,k;q)(-1)^k[r]_q^{m-k}[x+m]_q^n.$$
 \endproclaim
The generalized Euler numbers and polynomials of N$\ddot{o}$rlund's
type are defined by
$$\frac{2^r}{(e^{w_1t}+1)(e^{w_2t}+1)\cdots(e^{w_rt}+1)}=\sum_{n=0}^{\infty}E_n^{(r)}(x|w_1,
\cdots, w_r)\frac{t^n}{n!}, \tag33$$ and
$$E_{n}^{(r)}(w_1, \cdots, w_r)=E_{n}^{(r)}(0|w_1, \cdots, w_r).$$
Now, we can also consider  the $q$-extension of (33) as follows. For
$w_1, \cdots, w_r \in \Bbb Z_p$, and $\delta_1, \cdots, \delta_r
\in\Bbb Z, $we define
$$\aligned
&E_{n,q}^{(r)}(x|w_1, \cdots, w_r; \delta_1, \cdots, \delta_r)\\
&=\int_{\Bbb Z_p}\cdots\int_{\Bbb Z_p}[xw_1+\cdots+x_r w_r +x]_q^n
d\mu_{-q^{\delta_1}}(x_1)\cdots d\mu_{-q^{\delta_r}}(x_r),
\endaligned$$
and
$$E_{n,q}^{(r)}(w_1, \cdots, w_r;\delta_1, \cdots, \delta_r)=E_{n,q}^{(r)}(0|w_1, \cdots, w_r;\delta_1, \cdots,
\delta_r).$$ Thus, we have

$$E_{n,q}^{(r)}(x|w_1, \cdots, w_r;\delta_1, \cdots,
\delta_r)=\sum_{l=0}^n\binom{n}{l}(-1)^lq^{lx}\frac{(1+q^{\delta_1})\cdots
(1+q^{\delta_r})}{(1+q^{\delta_1+lw_1})\cdots (1+q^{\delta_r
+lw_r})}.\tag34$$ It seems to be interesting to consider another
$q$-extension of the N$\ddot{o}$rlund's type generalized Euler
numbers and polynomials as follows.
$$\aligned
&E_{n,q}^{*(r)}(x|w_1, \cdots, w_r;\delta_1, \cdots, \delta_r)\\
&=\int_{\Bbb Z_p}\cdots \int_{\Bbb
Z_p}[w_1x_1+\cdots+w_rx_r]_q^nq^{\delta_1
x_1+\cdots+\delta_rx_r}d\mu_{-1}(x_1)\cdots d\mu_{-1}(x_r),
\endaligned\tag35$$
and
$$E_{n,q}^{*(r)}(w_1, \cdots, w_r;\delta_1, \cdots, \delta_r)=E_{n,q}^{*(r)}(0|w_1, \cdots, w_r;\delta_1, \cdots,
\delta_r).$$
From (35), we note that
$$E_{n,q}^{*(r)}(x|w_1, \cdots, w_r;\delta_1, \cdots,
\delta_r)=2^r\sum_{l=0}^n\frac{\binom{n}{l}(-1)^lq^{lx}}{(1+q^{lw_1+\delta_1})\cdots(1+q^{lw_r+\delta_r})}.$$

\vskip 10pt

{\bf\centerline {\S 3. Further Remarks }} \vskip 10pt

For $h=0$, let us consider the polynomial of $E_{n,q}^{(0,r)}(x) $
and $E_{n, q}^{(0, -r)}(x)$ as follows.
$$E_{n,q}^{(0,r)}(x)=\int_{\Bbb Z_p}\cdots\int_{\Bbb
Z_p}[x_1+\cdots+x_r+x]_q^nq^{-\sum_{j=1}^rjx_j}d\mu_{-1}(x_1)\cdots
d\mu_{-1}(x_r),$$ and
$$E_{n,q}^{(0,-r)}=\sum_{l=0}^n
\frac {\binom{n}{l}(-1)^lq^{lx}}{\int_{\Bbb Z_p}\cdots \int_{\Bbb
Z_p}q^{l(x_1+\cdots+x_r)}q^{-\sum_{j=1}^rjx_j}d\mu_{-1}(x_1)\cdots
d\mu_{-1}(x_r)}.$$ Then, we have
$$\aligned
E_{n,q}^{(0,r)}(x)&=\frac{2^r}{(1-q)^n}\sum_{l=0}^n\frac{\binom{n}{l}(-1)^lq^{lx}}{(-q^{l-r};q)_r}\\
&=2^r\sum_{m=0}^{\infty}{\binom{m+r-1}{m}}_q(-1)^mq^{-rm}[m+x]_q^n,
\endaligned\tag36$$
and
$$\aligned
E_{n,q}^{(0,-r)}(x)&=\frac{1}{2^r
(1-q)^n}\sum_{l=0}^n\binom{n}{l}(-1)^lq^{lx}(-q^{l-r};q)_r\\
&=\frac{1}{2^r}\sum_{m=0}^r{\binom{r}{m}}_qq^{\binom{m}{2}}q^{-rm}[m+x]_q^n.
\endaligned\tag37$$
Let
$F_q^{(0,r)}(t,x)=\sum_{n=0}^{\infty}E_{n,q}^{(0,r)}(x)\frac{t^n}{n!}.$
Then we have
$$F_q^{(0,r)}(t,x)=2^r\sum_{m=0}^{\infty}{\binom{m+r-1}{m}}_q(-1)^m
q^{-rm}e^{[m+x]_q t}, $$ and
$$\aligned
F_q^{(0,-r)}(t,
x)&=\sum_{m=0}^{\infty}E_{n,q}^{(0,-r)}\frac{t^n}{n!}\\
&=\frac{1}{2^r}\sum_{m=0}^r{\binom{r}{m}}_qq^{\binom{m}{2}}q^{-rm}e^{[m+x]_qt}.
\endaligned$$
Let us consider the following polynomials.
$$\aligned
E_{n,q}^{(h,1)}(x)&=\int_{\Bbb
Z_p}q^{x_1(h-1)}[x+x_1]_q^nd\mu_{-1}(x_1)\\
&=\frac{2}{(1-q)^n}\sum_{l=0}^n\frac{\binom{n}{l}(-1)^lq^{lx}}{1+q^{l+h-1}}.
\endaligned\tag38$$
By the simple calculation of fermionic $p$-adic invariant integral
on $\Bbb Z_p$, we see that
$$\aligned
&q^x\int_{\Bbb Z_p}[x+x_1]_q^n q^{x_1(h-1)} d\mu_{-1}(x_1)\\
&=(q-1)\int_{\Bbb
Z_p}[x+x_1]_q^{n+1}q^{x_1(h-2)}d\mu_{-1}(x_1)+\int_{\Bbb
Z_p}[x+x_1]_q^nq^{x_1(h-2)}d\mu_{-1}(x_1).
\endaligned\tag39$$
From (39), we note that
$$q^xE_{n,q}^{(h,1)}(x)=(q-1)E_{n+1,q}^{(h-1,1)}(x)+E_{n,q}^{(h-1,1)}(x).$$

It is not difficult to show that
$$\int_{\Bbb
Z_p}q^{(h-1)x_1}[x+x_1]_q^nd\mu_{-1}(x_1)=\sum_{j=0}^n\binom{n}{j}[x]_q^{n-j}q^{jx}\int_{\Bbb
Z_p}[x_1]_q^jq^{(h-1)x_1}d\mu_{-1}(x_1).\tag40$$ By (40), we see
that
$$E_{n,q}^{(h,1)}(x)=\sum_{j=0}^n\binom{n}{j}[x]_q^{n-j}q^{jx}E_{j,q}^{(h,1)}=\left(q^xE_q^{(h,1)}+[x]_q \right)^n
,\text{ $n\geq 0$,}$$ where we use the technique method notation by
replacing $(E_q^{(h,1)})^n$ by $E_{n,q}^{(h,1)},$ symbolically. From
(25), we can also derive
$$\int_{\Bbb
Z_p}[x+x_1+1]_q^nq^{(x_1+1)(h-1)}d\mu_{-1}(x_1)+\int_{\Bbb
Z_p}[x+x_1]_q^nq^{(h-1)x_1}d\mu_{-1}(x_1)=2[x]_q^n.\tag41$$ Thus, we
see that
$$q^{h-1}E_{n,q}^{(h,1)}(x+1)+E_{n,q}^{(h,1)}(x)=2[x]_q^n.$$
For $x=0$, we have
$$q^{h-1}\left(qE_q^{(h,1)}+1\right)^n+E_{n,q}^{(h,1)}=2\delta_{n,0},$$
where $\delta_{n,0}$ is the Kronecker symbol.

It is easy to see that
$$E_{0,q}^{(h,1)}=\int_{\Bbb
Z_p}q^{x_1(h-1)}d\mu_{-1}(x_1)=\frac{2}{1+q^{h-1}}=\frac{2}{[2]_{q^{h-1}}}.$$
From (38), we note that
$$\aligned
E_{n, q^{-1}}^{(h,1)}(1-x)&=\int_{\Bbb Z_p}[1-x+x_1]_{q^{-1}}^n
q^{-x_1(h-1)}d\mu_{-1}(x_1)\\
&=(-1)^nq^{n+h-1}\frac{2}{(1-q)^n}\sum_{l=0}^n
\frac{\binom{n}{l}(-1)^lq^{lx}}{1+q^{l+h-1}}\\
&=(-1)^nq^{n+h-1}E_{n,q}^{(h,1)}(x).
\endaligned$$

In particular, for $x=1$, we have
$$E_{n,q^{-1}}^{(h,1)}(0)=(-1)^nq^{n+h-1}E_{n,q}^{(h,1)}(1)=(-1)^{n-1}q^nE_{n,q}^{(h,1)},
\text{ for $n\geq 1$.}$$

\vskip 20pt

\quad Taekyun Kim

\quad Division of General Education-Mathematics, Kwangwoon
University, Seoul 139-701, S. Korea \quad e-mail:\text{
tkkim$\@$kw.ac.kr}

\enddocument